\documentclass[12pt]{amsart}
\usepackage{epsfig}
\usepackage{amsmath}
\usepackage{amsfonts}
\usepackage{amssymb}
\usepackage{epigraph, fancyhdr}

\def\D{\mathcal D}
\def\F{\mathcal F}
\def\K{\mathcal K}
\def\J{\mathcal J}
\def\L{\mathcal L}
\def\H{\mathcal H}
\def\O{\mathcal O}

\def\1{\mathbf 1}
\def\M{{\mathcal M}}
\def\QQ{\mathbb Q}
\def\ZZ{\mathbb Z}
\def\CC{\mathbb C}

\def\Res{\operatorname{Res}}

\def\hat{\widehat}

\def\p{\partial}
\def\a{\alpha}
\def\b{\beta}

\def\gl{\lambda}

\def\lan{\langle}
\def\ran{\rangle}
\def\ev{\operatorname{ev}}

\def\td{\operatorname{td}}
\def\ch{\operatorname{ch}}
\def\qch{\operatorname{qch}}

\def\z{\zeta}

\renewcommand{\Delta}{\triangle}

\title[Explicit Reconstruction]
{Explicit Reconstruction in Quantum Cohomology and K-Theory}

\author[A. Givental]{Alexander GIVENTAL}
\thanks{This material is based upon work supported by the National 
Science Foundation under Grant DMS-1007164.} 
\dedicatory{To my friend Vadim Schechtman}

\date{June 25, 2014, revised May 12 and August 1, 2015}

\begin{document}



\begin{abstract}
Cohomological genus-0 Gromov-Witten invariants of 
a given target space can be encoded by the ``descendent potential,'' 
a generating function defined on the space of 
power series in one variable with coefficients in the cohomology space of the 
target. Replacing the coefficient space with the subspace multiplicatively 
generated by degree-2 classes, we explicitly reconstruct the graph of 
the differential of the restricted generating function from one point on it. 
Using the Quantum Hirzebruch--Riemann--Roch Theorem from our joint work 
\cite{GiT} with Valentin Tonita, we derive a similar reconstruction formula 
in genus-0 quantum K-theory. The results amplify the role in quantum cohomology
and quantum K-theory of the structures, based on divisor equations, of
$\D$-modules and $\D_q$-modules with respect to Novikov's variables. 
 
\end{abstract}

\maketitle


\section{Formulations}

Let $X$ be a compact K\"ahler (or, more generally, symplectic) manifold. Its 
genus-0 {\em descendent potential} is defined by
\[ \F (t) := \sum_{d \in \M} \sum_{n=0}^{\infty}\frac{Q^d}{n!}
\lan t(\psi), \dots, t(\psi)\ran_{0,n,d}, \]
where $\M \subset H_2(X,\ZZ)$ is the Mori cone of $X$, $Q^d$ stands for the
element corresponding to $d$ in the semigroup ring of $\M$,   
$t:=\sum_{k\geq 0} t_kz^k$ is a power series with coefficients $t_k$
which are cohomology classes of $X$,
and the correlator stands for the integral over the virtual fundamental class
$[X_{0,n,d}]$ of the moduli space of degree-$d$ stable maps to $X$ of rational 
curves with $n$ marked points:
\[ \lan \phi_1 \psi^{k_1},\dots, \phi_n\psi^{k_n}\ran_{0,n,d} :=
\int_{[X_{0,n,d}]} \ev_1^*(\phi_1) \psi_1^{k_1}\cdots \ev_n^*(\phi_n) \psi_n^{k_n}.\]
Here $\ev_i^*$ is the pull-back of cohomology classes from $X$ to $X_{0,n,d}$ by
the evaluation map at the $i$-th marked point, and $\psi_i$ is the 1st Chern
class of the line bundle over $X_{0,n,d}$ formed by cotangent lines to the curves at the $i$-th marked point.    
      
Following \cite{GiF}, we embed the graph of the differential of $\F$ into the 
{\em symplectic loop space} $\H$. By definition, it consists of formal 
$Q$-series whose coefficients are Laurent series in one 
indeterminate $z$ with vector coefficients from $H^*(X,\QQ)$.  

The ``loop space'' $\H$ (which is actually a $\ZZ_2$-graded module over the
Novikov ring $\QQ [[Q]]$) is equipped with the $\QQ [[Q]]$-valued even 
symplectic form
\[ \Omega (f,g):= \Res_{z=0} (f(-z),g(z))\ dz,\]
where $(\cdot, \cdot)$ is the Poincar{\'e} pairing (i.e.
$(a,b)=\int_X ab = \lan a, 1, b\ran_{0,3,0}$). 
Decomposing $\H$ into the sum $\H_{+}\oplus \H_{-}$ of complementary Lagrangian
subspaces (by the standard splitting of a Laurent series into the sum of 
the power $z$-series, and the polar part), we identify $\H$ with 
$T^*\H_{+}$. 
Translating the origin in $\H_{+}$ from $0$ to $-1z$ (the operation, refered to
as the {\em dilat{\'o}n shift}), we embed the graph of $d\F$ into 
$\H$ as a Lagrangian submanifold. Explicitly (see \cite{GiF}) it is given byXS the following {\em J-function}: 
\[ \H_{+} \ni t \mapsto \J (t):= -z + t(z)+\sum_{n,d,\alpha} \frac{Q^d}{n!} 
\phi_{\alpha} \lan \frac{\phi^{\alpha}}{-z-\psi}, t(\psi),\dots, t(\psi)
\ran_{0,n,d},\]
where $\{ \phi_{\alpha} \}$ and $\{ \phi^{\alpha} \}$ are Poincar{\'e}-dual
bases in $H^*(X,\QQ)$. 
  
In fact, this construction leads to some (rather mild) divergence problem. 
To elucidate it, pick a graded basis $\{ \phi_{\a} \}$ in $H^*(X,\QQ)$, 
and assume that $\phi_0=1$, and 
$\phi_{\a}$ with $\a =1,\dots, r = \operatorname{rk} 
H^2(X)$ are integer degree-2 classes $p_1,\dots,p_r$ taking non-negative 
values $d_i:=p_i(d)$ on degrees $d\in \M\subset H_2(X)$ of holomorphic curves 
in $X$. Writing 
\[ t_k=\sum_{\a}t_{k,\a}\phi_{\a} = t_{k,0}1+t_{k,1}p_1+\cdots+t_{k,r}p_r+
\text{the rest of the sum},\]
one can show (on the basis of string and divisor equations), that each 
$Q^d$-term in $\J$ contains the factor $e^{t_{0,0}/z}$ (which, unless expanded 
in powers of $t_{0,0}$, does not fit the space of formal Laurent series in $z$),
and besides comes with the factor $e^{\sum_i d_it_{0,i}}$ (which is not defined 
over $\QQ$). Also, as it follows from {\em dilaton equation}, with respect to
the variable $t_{1,0}$, the series has
convergence radius $1$. It follows from dimensional considerations that the 
rest of each $Q^d$-term is a polynomial in 
$1/z$ and in (finitely many of) the coefficients of the power series $t(z)$.

There are several ways to handle the problems. In this paper, we will ignore
the convergence properties by interpreting the J-function (and other geometric 
generating objects) in the sense of formal geometry. That is, $t\mapsto \J(t)$ 
is considered as the germ at $-z$ of a formal series in the components of 
the vector variables $t_k$ with coefficients which belong to the symplectic 
loop space. 

We will take $\QQ [[Q_1,\dots,Q_r]]$ on the role of the Novikov ring, 
and represent $Q^d$ by the monomial $Q_1^{d_1}\dots Q_r^{d_r}$. By virtue of the 
$Q$-adic convergence, one can specialize formal variables 
$t_{k,\a}$ to their values in the Novikov ring, taken from its maximal ideal
(which is necessary indeed in the case of $t_{0,0}, t_{0,i}$, and $t_{1,0}$). 
One can also make formal changes of the variables $\{ t_{k,\a} \}$ with
coefficients in the Novikov ring.  
          
\medskip

{\bf Theorem 1.} {\em Let $\sum_d I_d Q^d$, where $I_d(z,z^{-1})$ are 
cohomology-valued Laurent $z$-series, represent a point on the graph of 
$d\F$ in $\H$, and let $\Phi_{\alpha}$ be polynomials in $p_1,\dots, p_r$
(or, more generally, power $z$-series, with coefficients polynomial in 
$p_1,...,p_r$). Then the family  
\[ I(\tau):=
\sum_d I_dQ^d \exp \left\{ \frac{1}{z}\sum_{\alpha} \tau_{\alpha}
\Phi_{\alpha}(p_1-zd_1,\dots,p_r-zd_r) \right\}\]
lies on the graph of $d\F$.

Furthermore, for arbitrary scalar power series $c_{\alpha}(z)=\sum_{k\geq 0} 
\tau_{\alpha, k} z^k$, the linear combination
$\sum_{\alpha} c_{\alpha}(z) z\p_{\tau_{\alpha}} I$ of the derivatives 
also lies on the graph.

Moreover, in the case when $p_1,\dots,p_r$ generate the entire cohomology 
algebra $H^*(X,\QQ)$, and $\Phi_{\a}$ represent a linear basis, such linear
combinations comprise the whole graph.}

\medskip

{\tt Example 1.} Take $X=\CC P^{n-1}$, $p$ to denote the hyperplane class 
(and hence $p^n=0$), and $\Phi_i=p^i, i=0,...,n-1$, for a basis in $H^*(X,\QQ)$. The ``small J-function'' 
\[ (-z) \sum_{d\geq 0} \frac{Q^d}{(p-z)^n(p-2z)^n\cdots (p-dz)^n}\]
is known (see, for instance, \cite{GiH}) to represent a point on the graph
of $d\F$. It follows from THeorem 1 that the whole graph is comprises by
\[ (-z) \sum_{d\geq 0} \frac{Q^d e^{(\tau_0+\tau_1(p-dz)+\cdots+\tau_{n-1}(p-dz)^{n-1})/z}
\sum_{i=0}^{n-1}c_i(z)(p-dz)^i}
{(p-z)^n(p-2z)^n\cdots (p-dz)^n},\]
when $c_i(z)$ run arbitrary power series. More explicitly, one equates 
the power $z$-series part of this formula to $-z+t(z)$:  
\[ \sum_{i=0}^{n-1} (\tau_i -z c_i(z)) p^i + (\text{$Q$-adically small terms}) 
= -z + \sum_{i=0}^{n-1} p^i\sum_{k=0}^{\infty} t_{k,i}z^k,\] and expresses
$\tau_i$ and all coefficients of the series $c_i$ (here $c_0(0)$ needs to lie 
in a formal neighborhood of $1$) in terms of $\{ t_{k,i} \}$. Substituting these
expressions back into the formula, one obtains (according to Theorem 1) 
the standard form of the J-function for $\CC P^{n-1}$. 

\medskip

In K-theoretic version of GW-theory of a compact K\"ahler manifold $X$,
 the genus-0 descendent potential $\F^K$
is defined by the same formula as its cohomological counterpart:
\[ \F^K(t)=\sum_{d\in \M} \sum_{n=0}^{\infty}\frac{Q^d}{n!}
\lan t(L,L^{-1}),\dots, t(L,L^{-1})\ran^K_{0,n,d},\]
using the correlators
\[ \lan \Phi_1 L^{k_1}, \dots, \Phi_k L^{k_n}\ran^K_{0,n,d} :=
\chi (X_{0,n,d}; \O^{virt}\otimes \ev_1^*(\Phi_1)L_1^{k_1}\cdots 
\ev_n^*(\Phi_n)L_n^{k_n}).\]
Here $\chi$ is the holomorphic Euler characteristic (on $X_{0,n,d}$), 
$\O^{virt}$ is the virtual structure sheaf introduced by Yuan-Pin Lee 
\cite{YPLee}, $\Phi_i \in K^0(X)$ is a holomorphic vector bundle on $X$, 
$L_i^{k_i}$, $k_i\in \ZZ$, is the $k_i$th tensor power of the line bundle formed
by the cotangent lines to the curves at the $i$th marked point.
The input $t$ in $\F^K$ is a Laurent polynomial of $L$ with coefficients
in the K-ring of $X$. 

Adapting the symplectic loop space formalism, we embed the graph of $d\F^K$
as a Lagrangian submanifold into the ``space'' $\K$ consisting of power 
$Q$-series whose coefficients are rational functions in one indeterminate, $q$,
which take vector values in $K^0(X) \otimes \QQ$. 
Each rational function of $q$ is uniquely written as the sum of 
a Laurent polynomial and a rational function having no pole at $q=0$ and
vanishing at $q=\infty$. The space $\K$ is thereby decomposed into the direct 
sum of two subspaces, $\K_{+}$ and $\K_{-}$ respectively. They are Lagrangian
with respect to the symplectic form
\[ \Omega^K (f, g) =  \left[ \Res_{q=0}+\Res_{q=\infty} \right]\ (f(q), g(q^{-1})^K \frac{dq}{q},\]
where $(\cdot , \cdot)^K$ stands for the K-theoretic Poincar{\'e} pairing:
\[ (A,B)^K = \chi (X; A\otimes B) = \int_X \ch(A)\ch(B)\td(T_X).\]
Using this Lagrangian polarization to identify $\K$ with $T^*\K_{+}$, and 
applying the dilaton shift $1-q$, we identify the graph of $d\F^K$ with
a submanifold in $\K$, which is described explicitlyby the {\em J-function}
\begin{align*} &\K_{+}\ni t \mapsto \J^K(t)=\\
  &1-q+t(q,q^{-1})+\sum_{n,d,\alpha} \frac{Q^d}{n!}\Phi_{\alpha}
\lan \frac{\Phi^{\alpha}}{1-qL}, t(L,L^{-1}),\dots,t(L,L^{-1})\ran^K_{0,n+1,d}.\end{align*}
Here $\Phi_{\alpha}$ and $\Phi^{\alpha}$ run Poincar{\'e}-dual bases of $K^0(X)$.
Similar to the cohomological case, we consider $\J^K$ 
as a germ (at $1-q$) of a formal section of $T^*\H_{+}$. That is, it is a formal
series of the coordinates $t_{k,\a}$ (on the space of vector Laurent 
polynomials $\sum_{k,\a} t_{k,\a}\Phi_{\a} q^k$), whose coefficients are
$Q$-series with coefficients in rational functions of $q$.   

Let $P_1,...,P_r$ be line bundles over $X$ such that $c_1(P_i)=-p_i$, i.e. 
$d_i = - \int_d c_1(P_i)$.     
 
\medskip

{\bf Theorem 2.} {\em Let $\sum_d I_d Q^d$ be a point in $\K$, lying 
on the graph of $d\F^K$, and let $\Psi_{\a}$ be polynomials in $P_1,\dots, P_r$
(with coefficients which could be Laurent polynomials in $q$). Then the family 
\[ I^K(\tau)=\sum_d I_d Q^d \exp \{ \frac{1}{1-q}
\sum_{\alpha} \tau_{\alpha} \Psi_{\alpha}
(P_1q^{d_1}, \dots, P_rq^{d_r}) \}\]     
also lies on the graph. 

Furthermore, for arbitrary scalar Laurent polynomials $c_{\alpha}(q,q^{-1})$, 
the linear combinations $\sum_{\alpha} c_{\alpha}(q,q^{-1})(1-q) \p_{\tau_{\alpha}} I^K$ of the 
derivatives also lie on the graph.

Moreover, in the case when $P_1,\dots, P_r$ generate the algebra 
$K^0(X)\otimes \QQ$, and $\Phi_{\a}$ form a linear basis in it, such linear
combinations comprise the whole graph.}
   
\medskip 

{\tt Example 2.} Let $X$ be $\CC P^{n-1}$, $P = \O (-1)$ (and hence $(1-P)^n=0$), 
and $1, 1-P,\dots, (1-P)^{n-1}$ be the basis in $K^0(X)$. It was shown in 
\cite{GiL} that the following series lies on the graph of $d\F^K$:
\[ (1-q) \sum_{d=0}^{\infty}\frac{Q^d}{(1-Pq)^n(1-Pq^2)^n\cdots (1-Pq^d)^n}.\]
It follows that the whole graph can be parameterized this way:
\[ (1-q)\sum_{d=0}^{\infty}\frac{Q^d\ e^{\sum_{i=0}^{n-1}\tau_i (1-Pq^d)^i)/(1-q)}
\sum_{i=0}^{n-1} c_i(q,q^{-1}) (1-Pq^d)^i.}{(1-Pq)^n(1-Pq^2)^n\cdots(1-Pq^d)^n}.\]
More explicitly, one equates the Laurent polynomial part of this formula to
$(1-q)+t(q,q^{-1})$:
\[ \sum_{i=0}^{n-1}(1-P)^i(\tau_i+(1-q)c_i(q,q^{-1}))+
{\mathcal O}(Q)=1-q+\sum_{k,i}t_{k,i}q^k(1-P)^i\]
to express $\tau_i$ and all coefficients of the Laurent polynomials $c_i$
in terms of the variables $\{ t_{k,i} \}$. Substituting these expressions back
into the formula, one obtains the K-theoretic J-function of $\CC P^{n-1}$.     
 
{\tt Remark.} For target spaces, whose 2nd cohomology multiplicatively generate
the entire cohomology algebra, their cohomological and K-theoretic genus-0 
GW-invariants are reconstructible from small degree data, as it is established
by the reconstruction results of Kontsevich--Manin \cite{KM}, Lee--Pandharipande
\cite{LP}, and Iritani--Milanov--Tonita \cite{IMT}. Our results are closely 
related to them, and in a sense, explicize the reconstruction procedure.

{\tt Added at revision.} A preliminary version of this paper was posted
in May, 2014 to the author's website. After that we learned that Theorem 1
appears: (a) in a toric context --- in the earlier preprint \cite{CF-K},
Section 5.4, by I. Ciocan-Fontanine and B. Kim, and (b) in a general context
(though in a slightly less explicit form) --- in the even earlier
paper \cite{I}, Example 4.14, by H. Iritani. We are thankful to
Ionut Ciocan-Fontanine for this communication. 

\section{Proof of Theorem 1} 

Denote by $\L\subset \H$ the dilaton-shifted graph of $d\F$. 

\medskip

{\tt Step 1.} 
We begin by noting that modulo Novikov's variables, the graph is known
to have the form \cite{GiF}
\[ e^{-\tau/z} z\H_{+},\]\
where $\tau = \sum_{\a} \tau_{\a} \phi_{\a}$ runs the cohomology space of $X$. 

\medskip

{\tt Step 2.} The actual graph $\L$ is known (see Appendix 2 in \cite{CGi})
to have the form
\[ S_{\tau}^{-1}(z) z\H_{+},\]
where $\tau \mapsto S_{\tau}(z)$ is a certain family of matrices
(whose entries also depend on Novikov's variables), which has the following 
properties. Firstly, it is an $1/z$-series: $S=I+\O (1/z)$. 
Secondly, it belongs
to the ``twisted'' loop group: $S^{-1}(z)=S^*(-z)$, where ``$*$'' denotes 
transposition with respect to the Poincar{\'e} pairing. Thirdly, it is a 
fundamental solution to Dubrovin's connection on the tangent bundle
of the cohomology space of $X$:
\[ z\p_{\a} S = \phi_{\a}\bullet S,\]
where $\p_{\a}:=\p / \p \tau_{\a}$, and $\phi_{\a} \bullet$ is the matrix of 
quantum multiplication by $\phi_{\a}$ (it depends on the application point 
$\tau$ and on Novikov's variables, but not on $z$, and is self-adjoint).
Finally $S$ is constrained by the {\em string} and 
{\em divisor} equations. Namely, assuming as before, that $\{ \phi_{\a}\}$ 
is a graded basis in cohomology, with $\phi_0=1$ and $\phi_1=p_1,\dots, 
\phi_r=p_r$, we have:
\[ z\p_0 S = S, \ \ \text{and}\ \ z\p_i S = z Q_i \p_{Q_i} S + S p_i,\ \  
i=1,\dots, r.\]
(Here $p_i$ means the operator of multiplication by $p_i$ in the classical 
cohomology algebra of $X$.) 

Moreover, according to the ``descendent--ancestor correspondence'' theorem \cite{CGi} 
$S_{\tau}\L$ is tangent to $\H_{+}$ along $z\H_{+}$. This shows that
$\L$ is an {\em overruled Lagrangian cone}. By definition, this means that
tangent spaces $T_{\tau}$ to $\L$ (which are $S_{\tau}^{-1}\H_{+}$) are tangent
to $\L$ exactly along $zT_{\tau}$. We refer to \cite{CGi, GiF} for a more detailed discussion of this notion.

\medskip

{\tt Step 3.} Let $\D$ be the algebra of differential operators in 
Novikov's variables. It follows from the above divisor equations for $S$
that {\em tangent spaces $T_{\tau}=S_{\tau}^{-1}\H_{+}$ to $\L$ are $\D$-modules}
with respect to the action of $\D$ defined by the multiplication operators 
$Q_j$ and differentiation operators $zQ_i\p_{Q_i}-p_i$, where $p_i$ stands for 
multiplication by $p_i$ in the classical cohomology algebra of $X$. 
Consequently, the same is true about the {\em ruling spaces} $zT_{\tau}$.

Indeed, since $S^{-1}(-z)=S^*(z)$, and $p_i^*=p_i$, we have:
\[ (p_i+zQ_i\p_{Q_i}) S^{-1}(-z) = (Sp_i+zQ_i\p_{Q_i} S)^* = z \p_i S^* =
z\p_i S^{-1}(-z).\]
Here $S=S_{\tau}$ depends on $\tau \in H^*(X)$ and inependently on $Q$. Now
fix a value of $\tau=\tau(Q)$, and consider $f\in zT_{\tau} \subset \L$,
that is: $f = S^{-1}_{\tau} h$, where $h\in z\H_{+}$. 
Then
\[ (zQ_i\p_{Q_i}-p_i) f =
  z\p_iS_{\tau}^{-1}h + \sum_{\a} (Q_i\p_{Q_i}\tau_{\a}) z\p_{\a} S^{-1}_{\tau} h  +
  zS_{\tau}^{-1} (Q_i \p_{Q_i} h),\]
 where each summand on the right side lies in $zT_{\tau}$.

We arrive at the following conclusion.\footnote{This is a variant of
Lemma from the
proof of Quantum Lefschetz Theorem in Section 8 of \cite{CGi}. As we've
recently realized, the proof of it given in \cite{CGi} was incorrect.
Apparently, the argument was first corrected in \cite{CCIT} within the proof
of the orbifold version of the Quantum Lefschetz Theorem.}

\medskip

{\bf Lemma.} {\em Let $\Phi$ be a polynomial expression in $zQ_i\p_{Q_i}-p_i$.
Then the flow $f\mapsto e^{\epsilon \Phi/z}f$ preserves $\L$.}

\medskip

{\tt Proof.} If $T$ denotes the tangent space to $\L$ at $f\in zT$,
then $\Phi f/z \in T$, i.e. the linear vector field on $\H$: $f\mapsto \Phi f/z$
is tangent to $\L$.     

\medskip

{\tt Remark.} Since we are using differentiations in $Q$, it is 
counter-intuitive to think of Novikov's variables as constants. In fact one can 
think of the symplectic loop space $\H$ geometrically as the space of formal 
sections, over the spectrum of the Novikov ring, of the bundle whose 
fiber consists of Laurent $z$-series with vector coefficients. Likewise,
the cone $\L \subset \H$ consists of sections of the fibration whose fibers
are overruled Lagrangian cones. The differential operators $\Phi/z$ and their 
flows $e^{\epsilon \Phi/z}$ act by linear transformations on the space of sections
$\H$. In particular, $g=e^{\epsilon \Phi/z} f$ is an $\epsilon$-family of sections 
$Q\mapsto g(\epsilon, Q)$ of the fibration of overruled Lagrangian cones. 
One can choose any function $Q\mapsto \epsilon (Q)$ and obtain section 
$Q\mapsto g(\epsilon (Q), Q)$ lying in $\L$. We should note that it differs 
from $e^{\epsilon (Q) \Phi /z}f$ since multiplication by $\epsilon (Q)$ 
and $\Phi$ do not commute.    
  
\medskip

{\tt Step 4.} Write $f = \sum_d f_d Q^d$. Then 
\[ e^{\epsilon \Phi (\dots, p_i-zQ_i\p_{Q_p},\dots)/z}f= \sum_d f_d Q^d
e^{\epsilon \Phi (\dots, p_i-zd_i, \dots)/z},\]
which according to Step 3 lies in $\L$ whenever $f$ does. 
Here one can consider $\epsilon$ as a parameter, or take its value
from the Novikov ring (or at least from its maximal ideal). 

One obtains the first statement of Theorem 1 by replacing $\epsilon \Phi$
with a linear combination $\sum \tau_{\a} \Phi_{\a}$ of commuting differential
operators. 

The derivatives $\p_{\a} I(\tau)$ lie in the tangent space $T$ 
to $\L$ at $I_{\tau}$, and hence all linear combination 
$\sum c_{\a}(z) z\p_{\a} I(\tau)$, where $c_{\a}$ are scalar power $z$-series, 
lie in the same
ruling space $zT \subset \L$. 

When $p_1,\dots, p_r$ generate the entire cohomology algebra of $X$,
it follows from Step 1 that {\em modulo Novikov's variables}, 
such linear combinations comprise the whole of $\L$. Now the last statement of Theorem 1 follows the formal Implicit
Function Theorem.
   
\section{Proof of Theorem 2}

Let $\L^K \subset \K$ denote the graph of $d\F^K$. 

It is known (as explained in \cite{GiT}, Section 3) that $\L^K$ is 
an overruled Lagrangian cone too. 
More precisely, as in the case of quantum cohomology theory, there is 
a family $\tau \mapsto S_{\tau}(q, Q)$ of matrices depending on $\tau\in K^*(X)$
which transform $\L^K$ to $S_{\tau}\L^K$ tangent to $\K_{+}$ along
$(1-q)\K_{+}$. As a consequence, $\L^K$ is a cone whose tangent spaces 
$T_{\tau}=S_{\tau}^{-1}\K_{+}$ are $\QQ[q,q^{-1}]$ modules, and are 
tangent to $\L^K$ exactly along $(1-q)T_{\tau}$.
Theorem 2 is based on the property of the tangent and ruling spaces of $\L^K$
to be $D_q$-modules. Let us recall, following \cite{GiT}, how this is proved.
Another approach to this result is contained in \cite{IMT}. 

\medskip

The main result of \cite{GiT} (together with \cite{ToK,ToT}) is the Quantum Hirzebruch--Riemann--Roch Theorem
which completely characterizes $\L^K$ in terms of $\L$ in the following, 
``adelic'' way. For each complex value $\z$ of $q\neq 0$, one introduces the 
localization space $\K^{\z}$ which consists of series in $Q$ whose coefficients
are vectors in $K^0(X)\otimes \QQ(\z)$ and formal Laurent series in $1-q\z$. 
The {\em adelic map}
\[  \K  \to \hat{\K} := \prod_{\z} \K^{\z} \]
assigns to a rational function $f$ of $q$ the collection 
$(f^{(\z)})$ of its Laurent series expansions (one at each $q=\z^{-1}$).
 
Next, in each $\K^{\z}$, a certain cone $\L^{\z}$ is described. For $\z$ which
is not a root of $1$, $\L^{\z}=\K^{\z}_{+}$, the space of power series in $1-q\z$.
For $\z=1$, $\L^1 \subset \K^1$ is the graph of the differential of $\F^{fake}$,
the genus-0 descendent potential of {\em fake} quantum K-theory 
(studied in \cite{GiF, CGL, Co}). For $\z\neq 1$, which is a primitive $m$th 
root of $1$, $\L^{\z}$ is a certain linear subspace originating in a certain
fake twisted quantum K-theory with the orbifold target space $X/\ZZ_m$ (see
\cite{GiT} for more detail). 

The adelic characterization of $\L^K$ says that {\em $f \in \L^K$ if and only
if $f^{(\z)} \in \L^{\z}$ for each $\z$.}

Furthermore, $\L^{\z}$ have the following description in terms of the cone
$\L \subset \H$ of quantum cohomology theory. 

First, the {\em quantum Chern 
character} defines an isomorphism $\qch: \K^1\to \H^{even}$. By definition, $\qch$ 
acts by the usual Chern character on the coefficients of Laurent $q-1$-series, 
preserves Novikov's variables, and transforms $q$ into $e^z$. 
According to the ``quantum HRR theorem'' in {\em fake} quantum K-theory 
\cite{Co, CGL}, 
\[ \L^1 = \qch^{-1} \Delta \L, \ \ \text{where} \ \ 
\Delta \sim \prod_{\text{Chern roots $x$ of $T_X$}} \ 
\prod_{r=1}^{\infty} \frac{x-rz}{1-e^{-x+rz}}.\]
Here $\sim$ means taking the ``Euler-Maclaurin asymptotic'' of the right hand side. 
We won't remind the reader what it is (see, for instance, \cite{CGi, GiT}), but
note that (just as the expression on the right hand side suggests), it is 
multiplication by a series in $z^{\pm 1}$ built of operators of multiplication
in classical cohomology algebra of $X$, but {\em 
independent of Novikov's variables}. As a consequence, all tangent and ruling
spaces of the overruled Lagrangian cone $\L^{fake}$ are $\D$-modules (just
like those of $\L$ are). 

Then, when $\z\neq 1$ is a primitive $m$th root of $1$, one can give the 
following (somewhat clumsy) description of $\L^{\z}$.
On the cone $\L^K$, there is the point, denoted 
$\J(0)$ which corresponds to the input $t=0$. It is called the ``small 
J-function,'' and modulo $\K_{-}$, it is congruent to the dilaton shift $1-q$. 
Expanding $\J(0)$ into a Laurent series in $q-1$, we obtain the corresponding
point $\J(0)^{(1)}$ in $\L^{fake}$. 
The tangent space to $\L^{fake}$ at $\J(0)^{(1)}$ has the form 
\[   \Delta (z) S_{\tau(Q)}^{-1}(z, Q) \H_{+}^{even},\ \text{where} \ z=\log q,\]
and $S_{\tau}$ is the $S$-matrix of the cohomological theory computed at
a certain value $\tau = \tau(Q)$ (characterized by the application point
$\J(0)^{(1)}$). In this notation, $f \in \L^{\z}$ if and only if for some $h\in \H_{+}^{even}$
\[ \qch(\nabla_{\z}\, f(q^{1/m}/\zeta)) =  \left( \Psi^m\, \Delta(m z) S_{\tau(Q^m)}^{-1}(m z, Q^m) \Psi^{1/m} \right) h, \]
where
\[ \nabla_{\z} =
e^{\textstyle -\sum_{k>0}\left( \frac{\Psi^k(T^*_X)}{k(1-\z^{-k}q^{k/m})} -
  \frac{\Psi^{km}(T^*_X)}{k(1-q^{km})} \right)},\] 
%
and $\Psi^k$ are the Adams operations $K^0(X)\to \K^0(X)$ acting (by way of
the Chern isomorphism) on cohomology classes of degree $2r$ as
multiplication by $k^r$.

Let $P_1,\dots, P_r$ be line bundles on $X$ such that $c_1(P_i)=-p_i$,
and let $\D_q$ be the algebra of finite-difference operators in Novikov's 
variables. By definition, it acts on $\K$ by the ``translation'' operators 
$P_iq^{Q_i\p_{Q_i}} = \exp (-p_i+(\log q) Q_i\p_{Q_i})$ and multiplications by
$Q_j$. Let us show that $\L^{\z}$ is a $\D_q$-module.

Indeed, first note that that $\nabla_m$ and $\Delta$ commute with $\D$ (and hence with $\D_q$). Next, on functions of $Q^m$, we have
\[ \z^{-Q_i\p_{Q_i}}=\z^{-m Q_i^m\p_{Q_i^m}}=(\z^{-m})^{Q_i^m\p_{Q_i^m}}=1.\]
Therefore
\[(q^{1/m}/\z)^{Q_i\p_{Q_i}}e^{-p_i}=q^{Q_i^m\p_{Q_i^m}}e^{-p_i}=e^{-p_i+zQ_i^m\p_{Q_i^m}}.\]
Since $p_i\Psi^m=\Psi^mp_i/m$, the operator $zQ_i^m\p_{Q_i^m}-p_i$, when commuted accross $\Psi^m$, becomes $zQ_i^m\p_{Q_i^m}-p_i/m$. By the divisor equations for $S$, we have: 
\[ (z Q_i^m \p_{Q_i^m}-p_i/m) S_{\tau}^{-1}(m z, Q^m) =
z\p_i S_{\tau}^{-1}(m z, Q^m).\]
The remaining part of the computation (such as differentiation in $Q$
hidden in $\tau=\tau(Q^m)$) works out the same way as in the cohomological case.
It is also essential here that $\H_{+}^{even}$ is invarian under any differential of finite difference opeartors, including $\z^{Q_i\p_{Q_i}}$. 

\medskip

 As it was explained in \cite{GiT}, the adelic characterization of $\L^K$ now
implies that all tangent spaces $T_{\tau}$ to $\L^K$ (as well
as the ruling subspaces $(1-q)T_{\tau} \subset \L^K$ of the cone) are 
$\D_q$-modules.

Indeed, the whole space $\K$ is $\D_q$-invariant. If $f\in \K$ lies in a ruling
space $(1-q)T_f \subset \L^K$, then the adelic components $f^{(\z)}$ lie
in $\L^{\z}$. By the previous discusssion, $P_iq^{Q_i\p_{Q_i}} f^{(\z)} \in \L^{\z}$.
Therefore, by the adelic characterization, $P_iq^{Q_i\p_{Q_i}} f \in (1-q)T_f$.

\medskip

The proof of Theorem 2 proceeds now the same way as that of Theorem 1. If
$\Psi (\dots, P_iq^{Q_i\p_{Q_i}},\dots)$ is a polynomial expression in the
translation operators, then the linear vector field on $\K$ given by
$f\mapsto \Psi f/(1-q)$ is tangent to the cone $\L^K$, and therefore
the flow $f \mapsto e^{\epsilon \Psi /(1-q)} f$ preserves $\L^K$.

Decomposing $f$ into $Q$-series $\sum_d f_d Q^d$, we find that 
\[ e^{\epsilon \Psi /(1-q)}f = \sum_d f_d Q^d e^{\epsilon \Psi (\dots, P_iq^{d_i},\dots)/(1-q)}.\]
Replacing $\epsilon \Psi$ 
with a linear combination $\sum \tau_{\a} \Psi_{\a}$ of finite difference 
operators, one concludes that the family $\tau \mapsto I^K(\tau)$ (introduced
in Theorem 2) lies in $\L^K$. The derivatives $\p_{\a} I(\tau)$ lie
in the tangent space $T$ to $\L^K$ at $I^K(\tau)$. Since
$T$ is a module over $\QQ [q,q^{-1}]$, and $(1-q)T\subset \L^K$, one
finds that $\sum c_{\a}(q,q^{-1})(1-q)\p_{\a} I^K(\tau)$ also lie in $\L^K$. 
Finally, assuming that $P_1,\dots, P_r$ 
generate $K^0(X)$, one derives that such linear combinations comprise the 
whole of $\L^K$ by checking this statement modulo Novikov's variables, and 
employing the formal Implicit Function Theorem. 
    
\section{Further implications and generalizations}

{\bf A. Birkhoff's factorization and mirror maps.} When $H^*(X,\QQ)$ is 
generated by the degree-2 classes $p_1,\dots, p_r$, Theorems 1 and 2 can be 
reformulated as the following reconstruction results for the ``S-matrix.''
Starting with polynomials $\Phi_{\a}(p)$ representing a basis in $H^*(X,\QQ)$, 
and with a point $\sum I_d Q^d$ on the cone $\L$, one obtains a family
of such points 
\[ I(\tau)=\sum I_d(z,z^{-1}) Q^d e^{-\sum \tau_{\a}\Phi_{\a}(p-dz)/z}.\]
We may assume here that $I_0=-z$. The derivatives $\p_{\a} I$ form a 
$\QQ[[z]]$-basis in the tangent spaces to 
$\L$ (depending on $\tau$). The square matrix $U:=
[(\p_{\a}I, \phi^{\b})]$, 
formed by the components of these derivatives, can be factored into the product
of $U(z,z^{-1})=V(z)W(z^{-1})$ of two matrix series (in the variables $\tau$ 
and $Q$), 
whose coefficients are power series of $z$ (on the left) and polynomial 
functions of $z^{-1}$ (on the right). In the procedure (known as Birkhoff 
factorization), one may assume that $W(0)=I$.  
Then $W$ coincides with $S_{\tau}(-z^{-1})$ up to a change of variables
$\tau_{\a} \mapsto \tau_{\a}+{\mathcal O}(Q)$, which generalizes the 
``mirror map'' known in the mirror theory. To describe the change of variables,
assume that $\Phi_0=1$, and note that the ``first row'' of $W$ 
has the form 
\[ 1-z^{-1}\sum \Phi_{\a}(p)(\tau_{\a}+{\mathcal O}(Q))+o(z^{-1}).\]
The mirror map is read off the $z^{-1}$-term of the expansion.

In quantum K-theory, a similar result is obtained by Birkhoff factorization
$U=VW$, where the entries of $U$, $V$, and $W$ are built respectively of 
arbitrary rational functions, Laurent polynomials, and reduced rational 
functions of $q$ regular at $q=0$.  
  
\medskip  
 
{\bf B. Torus-equivariant theory.} 
It is often useful \cite{GiG} to consider GW-invariants 
{\em equivariant} with respect to a torus action on $X$.
The above results apply to this case without any significant changes.
One only needs to extend the coefficient ring by the power series completion
$\QQ [[\lambda]]$ of the coefficient ring of the equivariant theory. 
For example, when the torus $T^n$ of diagonal matrices acts on 
$X=\CC P^{n-1}=\operatorname{proj}(\CC^n)$, the $T^n$-equivariant cohomology
algebra of $X$ is described by the relation $(p-\gl_1)\cdots (p-\gl_n)=0$. 
For the purpose of employing fixed-point localization, it
is convenient to assume that the hyperplane class $p$ localizes to each of 
the values $\gl_j$. However, for the purpose of our proof it suffices to 
assume that $\gl_j$ are generators of the formal series ring 
$\QQ[[\gl_1,\dots,\gl_n]]$, and obtain the following 
parameterization of the graph $d\F$ in the $T^n$-equivariant GW-theory:
\[ (-z) \sum_{d\geq 0} \frac{Q^d e^{(\tau_0+\tau_1(p-dz)+\cdots+\tau_{n-1}(p-dz)^{n-1})/z}
\sum_{i=0}^{n-1}c_i(z)(p-dz)^i}
{\prod_{j=1}^n (p-\gl_j-z)(p-\gl_j-2z)\cdots (p-\gl_j-dz)}.\]
Here the fractions $1/(p-\gl-rz)$ are to be interpreted as Laurent polynomials in $z^{-1}$ modulo high powers of $\gl$.

\medskip

{\bf C. Twisted GW-invariants.}
Our results also extend to the case of twisted GW-invariants in the sense of
\cite{CGi} (e.g. ``local'' ones, i.e.
GW-invariants of the non-compact total space of a vector bundle $E\to X$, 
or GW-invariants of the ``super-bundle'' $\Pi E \to X$, which in genus $0$ 
are closely related to those of the zero locus of a section of $E$).
In such cases, to remove degenerations caused by non-compactness, one needs 
to act equivariantly, equipping $E$ with the fiberwise scalar circle action. 
To adapt our arguments to this case, it suffices to work over the coefficient 
ring $H^*(BS^1,\QQ)=\QQ[\gl]$ localized to $\QQ ((\gl))$. For example,
the graph of $d\F$ of the local theory on the total space $E$ of 
degree $l$ line bundle over $\CC P^{n-1}$, for $l>0$
obtains the following description:
\[ (-z) \sum_{d\geq 0} 
\frac{Q^d e^{(\tau_0+\tau_1(p-dz)+\cdots+\tau_{n-1}(p-dz)^{n-1})/z}
\sum_{i=0}^{n-1}c_i(z)(p-dz)^i}{\prod_{r=0}^{ld}(lp+\gl-rz)\ 
\prod_{r=1}^d (p-rz)^n}. \]
Here $p^n=0$, while the fractions $1/(lp+\gl-rz)$ should be expanded as  
power $z$-series, whose coefficients, however, can be Laurent series of 
$\gl$.

\enddocument